\documentclass[11pt]{article}
\usepackage{amsmath}
\usepackage{amssymb}
\usepackage{bbm}
\usepackage[margin=2.5cm]{geometry}
\usepackage{hyperref}
\usepackage{tikz}
\usetikzlibrary{decorations.pathmorphing,arrows}
\usepackage{braket}
\usepackage{setspace}

\newtheorem{thm}{Theorem}[section]
\newtheorem{cor}[thm]{Corollary}
\newtheorem{prop}[thm]{Proposition}

\newtheorem{lemma}[thm]{Lemma}

\newenvironment{proofof}{}{\hfill$\square$\vskip.5cm}
\newcommand{\R}{\mathbb{R}}
\newcommand{\N}{\mathbb{N}}
\newcommand{\C}{\mathbb{C}}

\title{About the Hardy-Ramanujan partition function asymptotics
}
\author{
Shannon Starr\footnote{\href{mailto:slstarr@uab.edu}{slstarr@uab.edu}}
\\
\large Department of Mathematics\\
\large University of Alabama at Birmingham\\
\small 1402 10th Avenue South,
Birmingham, AL 35294--1241
}
\date{\today}

\onehalfspace

\begin{document}

\maketitle

\abstract{The Hardy-Ramanujan partition function asymptotics is a famous result in the asymptotics of combinatorial
sequences. It was originally derived using complex analysis and number-theoretic ideas by Hardy and Ramanujan.
It was later 
re-derived by Paul Erd\H{o}s using real analytic methods. Later still, D.J.E.~Newman used just the usual Hayman saddle-point approach,
ubiquitous in asymptotic analysis.
Fristedt introduced a probabilistic approach, which was further extended by Dan Romik, for restricted partition functions.
Our perspective is that the Laplace transform changes the essentially algebraic generating function into an exponential
form. Using this, we carry out the exercise of deriving the leading order asymptotics,  following the Fristedt-Romik approach.
We also give additional examples of the Laplace transform method. }

\section{Introduction}

During the course of exploring a moment method approach to the generalized Ulam problem, the author
found it useful to exploit the Laplace transform in the extraction of asymptotics from a generating function
related to intersections of a pair of random walks \cite{HosseinStarr}.
We then considered some additional examples of the method. In the next section we illustrate with Newton's binomial formula.
In the proceeding section, we consider the $q$-binomial formula. But our main example, which is suprisingly simpler in some ways, is the partition function.

In \cite{HardyRamanujan} Hardy-Ramanujan derived leading-order asymptotics (as well as a full asymptotic series) of the partition function
$$
p(n)\, =\, |\operatorname{Par}(n)|\, =\, |\{\nu : \{1,2,\dots\} \to \{0,1,\dots\}\, :\, \sum_{k=1}^{\infty} k \nu(k) = n\}|\, .
$$
The generating function, for $|z|<1$ is 
$$
\sum_{n=0}^{\infty} p(n) z^n\, =\, \frac{1}{(z;z)_{\infty}}\, =\, \frac{1}{\prod_{k=1}^{\infty} (1-z^k)}\, .
$$
One can change this from an algebraic form to an exponential form by using the Laplace transform
$$
\frac{1}{\prod_{k=1}^{\infty} (1-z^k)}\, =\, \int_{\R_+^{\N}} e^{-\sum_{k=1}^{\infty} t_k (1-z^k)}\, \prod_{k=1}^{\infty} dt_k\, .
$$
Doing this, the exercise of extracting the leading-order asymptotics follows essentially by the Hayman saddle-point method, in a manner which seems easier than 
that of D.J.E. Newman \cite{Newman}.
By using exponential random variables instead of writing out the full exponential integrals, a connection can be seen to Fristedt's method \cite{Fristedt},
which was carried further by Romik \cite{Romik}.
Essentially, this makes the example of the Hardy-Ramanujan leading order asymptotics an exercise.
But we carry out this exercise.

\subsection{Digression: Stirling's formula review}
We recall one of the derivations of Stirling's formula, which will be useful for us, later. Starting from
$$
\ln(n!)\, =\, \sum_{k=1}^{n} \ln(k)\, ,
$$
and doing the first step of the Euler-Maclaurin summation formula, we obtain
$$
\ln(n!)\, =\, \sum_{k=1}^{n} \int_{k-\frac{1}{2}}^{k+\frac{1}{2}} \ln(x)\, dx
- \sum_{k=1}^{n} \int_{k-\frac{1}{2}}^{k+\frac{1}{2}} (\ln(x)-\ln(k))\, dx\, .
$$
Noting that the second summation is a partial sum of a summable series, we note that the first sum is 
$$
\sum_{k=1}^{n} \left(\left(k+\frac{1}{2}\right) \ln\left(k+\frac{1}{2}\right) - \left(k-\frac{1}{2}\right) \ln\left(k-\frac{1}{2}\right) - 1\right)\, .
$$
This motivates considering
$$
\frac{n!}{(n+\frac{1}{2})^{n+\frac{1}{2}}}\, \left(\frac{1}{2}\right)^{1/2}\, e^n\,
=\, \exp\left(\sum_{k=1}^{n} \left(\ln(k) - \left(k+\frac{1}{2}\right) \ln\left(k+\frac{1}{2}\right) + \left(k-\frac{1}{2}\right) \ln\left(k-\frac{1}{2}\right) + 1\right)\right)
$$
But the right-hand-side is easily seen to be 
$$
\exp\left(-\frac{1}{2}\sum_{k=1}^{n} \ln\left(1-\frac{1}{4k^2}\right)\right)
\exp\left(\sum_{k=1}^{n}\left(1-k\, \ln\left(1+\frac{1}{2k}\right) + k\, \ln\left(1-\frac{1}{2k}\right) \right)\right)\, .
$$
Both of the sums are summable, and exactly calculable, using the product formula for sine:
$$
\sin(\pi x)\, =\, \pi x \prod_{k=1}^{n} \left(1- \frac{x^2}{k^2}\right)\, .
$$
Setting $x=1/2$ and taking logarithms, we see that
\begin{equation}
\label{eq:1}
-\sum_{k=1}^{\infty} \ln\left(1-\frac{1}{4k^2}\right)\, =\, \ln\left(\frac{\pi}{2}\right)\, .
\end{equation}
Moreover, the integral of $\ln(\sin(\pi x)/(\pi x))$ can be expressed in terms of the dilogarithm,
$$
\operatorname{Li}_2(z)\, =\, \sum_{k=1}^{\infty} \frac{z^k}{k^2}\, =\, -\int_0^z \frac{\ln(1-\xi)}{\xi}\, d\xi\, . 
$$ 
(See for example Zagier \cite{Zagier}.)
More precisely, using $\operatorname{Li}_2(1) = \zeta(2)=\pi^2/6$, we have
$$
	\int_0^x \ln\left(\frac{\sin(\pi x)}{\pi x}\right)\, dx\,
	=\, -\frac{i\pi}{12} + x \left(\frac{i \pi}{2} - \ln(2)\right) - \frac{i \pi}{2}\, x^2
- (x \ln(\pi x)-x) - \frac{1}{2\pi i}\, \operatorname{Li}_2(e^{2\pi i x})\, .
$$
Then, using the fact that $\operatorname{Li}_2(-1) = -\pi^2/12$, we may deduce
\begin{equation}
\label{eq:2}
\sum_{k=1}^{\infty} \frac{\ln(1+\frac{1}{2k})-\ln(1-\frac{1}{2k})-\frac{1}{k}}{1/k}\, =\,  \frac{1}{2}\, \left(1-\ln(2)\right)\, .
\end{equation}
Therefore, since $(n+\frac{1}{2})^{n+\frac{1}{2}} = n^{n+\frac{1}{2}} \left(1+\frac{1}{2n}\right)^{n+\frac{1}{2}} \sim e^{1/2} n^{n} \sqrt{n}$,
as $n \to \infty$,
we obtain Stirling's formula
$$
\lim_{n \to \infty} \frac{n!}{n^n e^{-n} \sqrt{2\pi n}}\, =\, 1\, .
$$
The reason for reviewing this derivation of Stirling's formula is that these equations will return.

\section{First warm-up example: Newton's binomial formula}
We start with a warm-up example, to demonstrate the usual idea of the Laplace transform.
The idea is to replace an algebraic generating function with a more exponential one, which works
better with the saddle point method.

Perhaps the most basic generating function is Newton's binomial formula. Suppose we have a number $s \in (0,\infty)$.
Then
$$
	\frac{1}{(1-z)^s}\, =\, \sum_{n=0}^{\infty} \frac{(s+n-1)_n}{n!}\, z^n\, ,
$$
for each $|z|<1$, where the falling factorial Pocchammer symbol is $(x)_n = x(x-1)\cdots(x-n+1)$.
Therefore, by Cauchy's integral formula, we know that
$$
	\frac{(s+n-1)_n}{n!}\, =\, \oint_{\mathcal{C}(0;r)} \frac{1}{z^n} \cdot \frac{1}{(1-z)^s}\, \cdot \frac{dz}{2\pi i z}\, ,
$$
for all $r \in (0,1)$. The optimal choice would be $r=n/(n+s)$. So
$$
	\frac{(s+n-1)_n}{n!}\, =\, \left(1+\frac{s}{n}\right)^{n} \int_{-\pi}^{\pi} e^{-in\theta} \cdot \frac{1}{\left(1-r e^{i\theta}\right)^s}\, \cdot \frac{d\theta}{2\pi}\, ,
$$
where we allow ourselves to continue to use $r$ in place of $n/(n+s)$, where convenient.
Then, we use the Laplace transform identity
$$
	\operatorname{Re}[a]>0\qquad \Rightarrow\qquad \frac{1}{a^s}\, =\, \int_0^{\infty} e^{-at}\, \frac{t^s}{\Gamma(s)}\, \frac{dt}{t}\, .
$$
This gives
$$
	\frac{(s+n-1)_n}{n!}\, =\, \frac{1}{\Gamma(s)}\, \left(1+\frac{s}{n}\right)^{n} \int_{-\pi}^{\pi} e^{-in\theta}
 \left(\int_0^{\infty} e^{-t(1-re^{i\theta})}\, t^s\, \frac{dt}{t}\right)
\, \frac{d\theta}{2\pi}\, .
$$
By the Fubini-Tonelli theorem, this may be rewritten as 
$$
	\frac{(s+n-1)_n}{n!}\, =\, \frac{1}{\Gamma(s)}\, \left(1+\frac{s}{n}\right)^{n} \int_0^{\infty} e^{-t}\, t^{s-1}
\left(\int_{-\pi}^{\pi} e^{-in\theta} e^{tre^{i\theta}}\, \frac{d\theta}{2\pi}\right)\, dt\, .
$$
Let us summarize some elementary calculations as follows:
\begin{lemma}
\label{lem:analytic}
For any $0<\epsilon<\delta<1$, if $t<n^{1-\delta}$, then
$$
\int_{-\pi}^{\pi} e^{-in\theta} e^{tre^{i\theta}}\, \frac{d\theta}{2\pi}\, =\, \frac{(tr)^n}{n!}\, =\, O(n^{-\epsilon n})\, ,\
\text{ as $n \to \infty$.}
$$
\end{lemma}
\begin{lemma}
\label{lem:0a}
Letting $\tau = t/(n+s)$, and using $e^{i\theta}-1 = i \sin(\theta) - 2 \sin^2(\theta/2)$, we have
\begin{equation*}
\begin{split}
&\hspace{-0.75cm}e^s \int_0^{\infty} e^{-t}\, t^{s-1}
\left(\int_{-\pi}^{\pi} e^{-in\theta} e^{tre^{i\theta}}\, \frac{d\theta}{2\pi}\right)\, dt\\
&\hspace{0.75cm}=\, (n+s)^s \int_0^{\infty} e^{-s (\tau-1)}\, \tau^{s-1}
\left(\int_{-\pi}^{\pi} e^{-in(\theta-\tau \sin(\theta))} e^{-2 n \tau \sin^2(\theta/2)}\, \frac{d\theta}{2\pi}\right)\, d\tau\, .
\end{split}
\end{equation*}
\end{lemma}
\begin{cor}
For $s>0$ fixed, and for any $0<\epsilon<\delta<1$, we have
$$
\frac{(s+n-1)_n}{n!} \cdot \frac{\Gamma(s)}{n^s}\, 
=\, O(n^{-\epsilon n + s})+
\int_{n^{-\delta}}^{\infty} e^{-s (\tau-1)}\, \tau^{s-1}
\left(\int_{-\pi}^{\pi} e^{-in(\theta-\tau \sin(\theta))} e^{-2 n \tau \sin^2(\theta/2)}\, \frac{d\theta}{2\pi}\right)\, d\tau\, ,
$$
as $n \to \infty$.
\end{cor}
Note that changing variables to $\theta = \Theta/\sqrt{n}$, for $\Theta = o(n^{1/4})$ and $\tau=O(1)$, we have
$$
e^{-in(\theta-\tau \sin(\theta))} e^{-2 n \tau \sin^2(\theta/2)}\, \sim\, e^{in^{1/2}(1-\tau)\Theta} e^{-\tau \Theta^2/2}\, .
$$
Integrating the right-hand-side with respect to $\Theta$ over all of $\R$ would gives $e^{-n(1-\tau)^2/(2\tau)}$
times $\sqrt{2\pi}$.
But there are some technical details. So we state the following as a lemma,  which we will prove in Appendix \ref{app:first}.
\begin{lemma}
\label{lem:Hay1}
By a variation of the usual Hayman saddle-point method, allowing for deforming the contour of integration into the complex plane,
for $\tau > n^{-\delta}$, and $\delta \in (0,1/3)$ we have
$$
\int_{-\pi}^{\pi} e^{-in(\theta-\tau \sin(\theta))} e^{-2 n \tau \sin^2(\theta/2)}\, \frac{d\theta}{2\pi}\, 
=\, \frac{e^{n(\ln(\tau)+1-\tau)}}{ \sqrt{2\pi n}}\, \left(1+o(1)\right) + o\left(\frac{1}{n}\right)\, ,\ \text{ as  $n \to \infty$.}
$$
\end{lemma}
This leads to the following lemma, which we will also prove in Appendix \ref{app:first}.
\begin{lemma}
\label{lem:Hay2}
By the usual Hayman saddle-point method, 
$$
\int_{n^{-\delta}}^{\infty} e^{-s(\tau-1)} \tau^{s-1}\, \frac{e^{n(\ln(\tau)+1-\tau)}}{\sqrt{2\pi n }}\, d\tau\,
\sim\, \frac{1}{n}\, ,\ \text{ as  $n \to \infty$.}
$$
\end{lemma}
Hence, we have derived the well-known result, following. (See, for instance, Definition 1.1.1 
and Corollary 1.1.5 in \cite{AndrewsAskeyRoy}.) We do provide a short proof in Appendix \ref{app:first}.
\begin{cor}
\label{cor:Gamma}
For $s>0$, we have
$$
\frac{(s+n-1)_n}{n!} \cdot \frac{\Gamma(s)}{n^s}\, \sim\, \frac{1}{n}\, ,\ \text{ as  $n\to\infty$.}
$$
In other words,
$$
\Gamma(s)\, =\, \lim_{n \to \infty} \frac{n!\, n^{s-1}}{(n+s-1)_n}\, .
$$
\end{cor}
We prove the last two lemmas in the appendix because they really do follow the usual Hayman analysis, wherein a saddle-point
is well-approximated by a Gaussian. (The corollary is shorter.)
What is somewhat novel is the use of the Laplace transform.

\section{Second warm up example: the $q$-Stirling formula}

By the $q$-binomial formula, we know that
$$
	\sum_{n=0}^{\infty} \frac{(a;q)_n}{(q;q)_n}\, z^n\, =\, \frac{(az;q)_{\infty}}{(z;q)_{\infty}}\, ,
$$
for $|z|<1$ and $|q|<1$. See, for example,Section 1.3 of Gasper and Rahman \cite{GasperRahman}.
The simplest case of this is to set $a=0$ and obtain the $q$-deformed version of the exponential
$$
	\sum_{n=0}^{\infty} \frac{1}{(q;q)_n}\, z^n\, =\, \frac{1}{(z;q)_{\infty}}\, .
$$
The asymptotics of the coefficients $1/(q;q)_n$ have been considered by Moak.
But the method was the Euler-Maclaurin formula. We also use the Euler-Maclaurin step.
But here we wish to demonstrate the Laplace transform trick. So we re-derive one of Moak's results,
noting that there are also numerous other notable results in \cite{Moak}.
We also note that one interesting application of the $q$-Stirling formula is in \cite{StarrWalters}.

Here, the $q$-Pochhammer symbols are 
$$
(a;q)_n\, =\, \prod_{k=0}^{n-1} (1-aq^k)\ \text{ for $n \in \{0,1,\dots\}$, }\ \text{ and }\ (a;q)_{\infty}\, =\, \prod_{k=0}^{\infty} (1-aq^k)\, .
$$
Let us assume that $|a|<1$ so that all of the factors are within a distance 1 from $1$ in $\C$.
Then, for $\epsilon \to 0^+$, we have
\begin{equation}
\label{eq:backToThis}
\ln\left((a;e^{-\epsilon})_{\infty}\right)\, =\, \sum_{k=0}^{\infty} \ln(1-ae^{-\epsilon k})\,
\sim\, \frac{1}{\epsilon}\, \int_0^{\infty} \ln(1-e^{-t} a)\, dt\, ,
\end{equation}
and we also have
$$
\int_0^{\infty} \ln(1-e^{-t}a)\, dt\, =\, -\operatorname{Li}_2(a)\, .
$$
However, in order to have a more precise estimate of the sum, let us use one step of the Euler-Maclaurin
method:
$$
\sum_{k=0}^{\infty} \ln(1-ae^{-\epsilon k})\, =\, \int_{-1/2}^{\infty} \ln(1-a e^{-\epsilon x})\, dx
- \sum_{k=0}^{\infty} \int_{k-\frac{1}{2}}^{k+\frac{1}{2}} 
\left(\ln(1-ae^{-\epsilon x}) - \ln(1-a e^{-\epsilon k})\right)\, dx\, .
$$
We note that, changing variables to $t=\epsilon x$ and $dx=dt/\epsilon$, we have
\begin{equation*}
\begin{split}
\int_{-1/2}^{\infty} \ln(1-a e^{-\epsilon x})\, dx\, 
&=\, -\epsilon^{-1} \operatorname{Li}_2(a)
+ \epsilon^{-1} \int_{-\epsilon/2}^{0} \ln(1-a e^{-t})\, dt\\
&=\, -\epsilon^{-1} \operatorname{Li}_2(a) + \frac{1}{2}\, \ln(1-a) + o(1)\, ,
\end{split}
\end{equation*}
as $\epsilon \to 0$.
Moreover, by Taylor expansion, it is easy to see that 
$$
\sum_{k=0}^{\infty} \int_{k-\frac{1}{2}}^{k+\frac{1}{2}} 
\left(\ln(1-ae^{-\epsilon x}) - \ln(1-a e^{-\epsilon k})\right)\, dx\,
\sim\, -\frac{\epsilon^2}{96}\, \int_0^{\infty} \frac{dx}{\sinh^2(\frac{x}{2} - \frac{1}{2}\ln(a))}\, dx\, ,
$$
which is $O(\epsilon^2)$ as $\epsilon \to 0^+$, as long as $|a|<1$ is fixed.
So this correction may be neglected.
Hence, returning to equation (\ref{eq:backToThis}), we have
\begin{equation}
\label{eq:BackToIt2}
(a;e^{-\epsilon})_{\infty}\, \sim\, e^{-\epsilon^{-1} \operatorname{Li}_2(a)}\, \sqrt{1-a}\, ,
\end{equation}
as $\epsilon \to 0^+$.

Now using Cauchy's integral formula, for a fixed $\beta>0$ if we consider $(q;q)_n$ for $q=q_n=e^{-\beta/n}$ we have
$$
\frac{1}{(q;q)_n}\Bigg|_{q=\exp(-\beta/n)}\,
	=\, \int_{-\pi}^{\pi} \frac{1}{r^n e^{in\theta}}\, \frac{1}{(r e^{i\theta};e^{-\beta/n})_{\infty}}\, \frac{d\theta}{2\pi}\, ,
$$
for any $r \in (0,1)$. But setting $\theta=0$ and optimizing the integrand by choosing $r$, we find
$$
r\, =\, 1-e^{-\beta}\, .
$$
(We do this by using $\frac{d}{da}\, \operatorname{Li}_2(a) = -\ln(1-a)/a$.)
Therefore, we do have $|r|<1$.

Now, we introduce the Laplace transform by using IID exponential-1 random variables $\mathsf{T}_0,\mathsf{T}_1,\dots$, so that
$$
\forall t \in \R\, ,\ \text{ we have }\
	\mathbf{P}(\mathsf{T}_k \leq t)\, =\, \begin{cases} 0 & \text{ if $t<0$,}\\
1-e^{-t} & \text{ if $t\geq 0$.}
\end{cases}
$$
Then we have, for any $a_k$ such that $|a_k|<1$, the formula:
$$
\mathbf{E}[e^{a_k \mathsf{T}_k}]\, =\, \int_0^{\infty} e^{a_k t} e^{-t}\, dt\, =\, \frac{1}{1-a_k}\, .
$$
So we may write the infinite product in the Cauchy integral formula instead as 
$$
\frac{1}{(e^{-\beta/n};e^{-\beta/n})_n}\,
	=\, \frac{1}{r^n} \int_{-\pi}^{\pi} e^{-in\theta} 
\mathbf{E}\left[ \exp\left(r\sum_{k=0}^{\infty} e^{-\beta k/n} \mathsf{T}_k e^{i\theta}\right)\right]\, 
\frac{d\theta}{2\pi}\, .
$$
Now we factor out the integrand at $\theta=0$ by adding-and-subtracting in the exponent:
$$
\frac{1}{(e^{-\beta/n};e^{-\beta/n})_n}\,
	=\, \frac{1}{r^n} \int_{-\pi}^{\pi} e^{-in\theta} 
\mathbf{E}\left[ \exp\left(r\sum_{k=0}^{\infty} e^{-\beta k/n} \mathsf{T}_k\right)
\exp\left(r\sum_{k=0}^{\infty} e^{-\beta k/n} \mathsf{T}_k (e^{i\theta}-1)\right)\right]\, 
\frac{d\theta}{2\pi}\, .
$$
Now we use the 
exponential-tilting formula for the exponential random variable for $a_k \in (-1,\infty)$:
\begin{equation*}
\begin{split}
\mathbf{E}[e^{a_k\mathsf{T_k}}f(\mathsf{T_k})]\, 
&=\, \int_0^{\infty} e^{-t} e^{a_kt} f(t)\, dt\\
&=\, \int_0^{\infty} e^{-(1-a_k)t} f(t)\, dt\\
&=\, \int_0^{\infty} e^{-s} f\left(\frac{s}{1-a_k}\right)\, \frac{ds}{1-a_k}\\
&=\, \frac{1}{1-a_k}\, \mathbf{E}\left[f\left(\frac{\mathsf{T}_k}{1-a_k}\right)\right]\, .
\end{split}
\end{equation*}
So we have
\begin{equation}
\label{eq:temporaryPeak}
\begin{split}
\frac{1}{(e^{-\beta/n};e^{-\beta/n})_n}\,
	&=\, \frac{1}{r^n}\, \cdot \frac{1}{(r;e^{-\beta/n})_{\infty}} \\
&\qquad \cdot \int_{-\pi}^{\pi} e^{-in\theta}\, 
\mathbf{E}\left[ 
\exp\left(\sum_{k=0}^{\infty} \frac{re^{-\beta k/n}}{1-re^{-\beta k/n}} \mathsf{T}_k (e^{i\theta}-1)\right)\right]\, 
\frac{d\theta}{2\pi}\, .
\end{split}
\end{equation}
Note that we have already considered the asymptotics of the prefactor in equation (\ref{eq:BackToIt2}). What remains is to 
consider the asymptotics of the integral.

Now we will use two properties of the random variables to obtain control of the integral: Cram\'er's
large deviation principle to obtain basic bounds, and the central limit theorem for more precise results.

\subsection{Rough bounds from Cram\'er's large deviation principle}

We start with the following. (See for example
Dembo and Zeitouni \cite{DemboZeitouni}.)
\begin{prop}[Cram\'er's theorem] The sequence $(\mathsf{T}_1+\dots+\mathsf{T}_n)/n$ satisfies a large deviation principle
with large deviation rate function $\Lambda^*$, whose formula is 
$$
\Lambda^*(x)\, =\, x-1-\ln(x)\, ,
$$
the Legendre transform of the  cumulant generating function $\Lambda(t) = \ln\left(\mathbf{E}[e^{t\mathsf{T}}]\right)=-\ln(1-t)$.
\end{prop}
Analysis similar to the previous Euler-Maclaurin step shows that
\begin{equation}
\label{eq:EM2}
\begin{split}
\sum_{k=0}^{\infty} \frac{r e^{-k\epsilon}}{1-re^{-k\epsilon}}\, 
&\sim\,
\epsilon^{-1} \int_0^{\infty} \frac{r e^{-x}}{1-r e^{-x}}\, dx + \frac{1}{2}\cdot \frac{r}{1-r}\\
&=\, -\epsilon^{-1} \ln(1-r) + \frac{1}{2}\cdot \frac{r}{1-r}\, ,
\end{split}
\end{equation}
as $\epsilon \to 0$. 
Here we just want rough lower bounds. So we note that 
$$
\forall k \in \{0,\dots,n\}\, ,\ \text{ we have }\ 
\frac{r e^{-\beta k/n}}{1-re^{-\beta k/n}}\, \geq\, \frac{r e^{-\beta}}{1-re^{-\beta}}\, \geq\, r e^{-\beta}\, .
$$
This implies
$$
\sum_{k=0}^{\infty} \frac{re^{-\beta k/n}}{1-re^{-\beta k/n}} \mathsf{T}_k\, \geq\, r e^{-\beta} \sum_{k=0}^{n-1} \mathsf{T}_k\, .
$$
Hence by the Cram\'er theorem,
$$
\mathbf{P}\left(\sum_{k=0}^{\infty} \frac{re^{-\beta k/n}}{1-re^{-\beta k/n}} \mathsf{T}_k\leq re^{-\beta}n(1-x)\right)\,
\leq\, e^{-n\Lambda^*(x)(1+o(1))}\, .
$$
Note that we also have
\begin{equation*}
\begin{split}
\left|\mathbf{E}\left[ 
\exp\left(\sum_{k=0}^{\infty} \frac{re^{-\beta k/n}}{1-re^{-\beta k/n}} \mathsf{T}_k (e^{i\theta}-1)\right)\right]\right|\,
&\leq\, \mathbf{E}\left[ 
\exp\left(\sum_{k=0}^{\infty} \frac{re^{-\beta k/n}}{1-re^{-\beta k/n}} \mathsf{T}_k \operatorname{Re}[e^{i\theta}-1]\right)\right]\\
&=\, \mathbf{E}\left[ 
\exp\left(-2\sum_{k=0}^{\infty} \frac{re^{-\beta k/n}}{1-re^{-\beta k/n}} \mathsf{T}_k 
\sin^2\left(\frac{\theta}{2}\right)\right)\right]
\end{split}
\end{equation*}
Therefore, using $\sin^2(\theta/2) \leq \theta^2/\pi^2$ for $|\theta|<\pi$, we have that for $|\theta| \in [\theta_0,\pi]$, it is the case that
\begin{equation*}
\begin{split}
\left|\mathbf{E}\left[ 
\exp\left(\sum_{k=0}^{\infty} \frac{re^{-\beta k/n}}{1-re^{-\beta k/n}} \mathsf{T}_k (e^{i\theta}-1)\right)\right]\right|\,
\leq\, e^{-2re^{-\beta}n(1-x)\sin^2(\theta_0/2)} + e^{-n\Lambda^*(x)(1+o(1))}\, .
\end{split}
\end{equation*}
We will choose $\theta_0 = R_n/\sqrt{n}$ for $R_n = K\, \sqrt{\ln(n)}$ and sufficiently large fixed $K$. Then taking $x=1/2$, for example, we can make this quantity $O(1/n^p)$
for any desired $p<\infty$, if we choose $K$ sufficiently large (with $K=O(\sqrt{p})$ as $p \to \infty$).

One could probably construct a tailor-made argument, bespoke, without using Cram\'er.

\section{More precise estimates from the central limit theorem}

Now, for $\theta = \Theta/\sqrt{n}$ with $|\Theta|\leq R_n$, we use the formula
\begin{equation*}
\begin{split}
\mathbf{E}\left[ 
\exp\left(\sum_{k=0}^{\infty} \frac{re^{-\beta k/n}}{1-re^{-\beta k/n}} \mathsf{T}_k (e^{i\theta}-1)\right)\right]\,
&=\,
\exp\left(\sum_{k=0}^{\infty} \frac{re^{-\beta k/n}}{1-re^{-\beta k/n}} (e^{i\theta}-1)\right)\\
&\qquad \cdot \mathbf{E}\left[ 
\exp\left(\sum_{k=0}^{\infty} \frac{re^{-\beta k/n}}{1-re^{-\beta k/n}} (\mathsf{T}_k-1) (e^{i\theta}-1)\right)\right]\, ,
\end{split}
\end{equation*}
where we have explicitly added-and-subtracted the mean of each random variable in the exponent.
Writing $e^{i\theta}-1 = i \sin(\theta) - 2 \sin^2(\theta/2)$, we have
\begin{equation*}
\begin{split}
\exp\left(\sum_{k=0}^{\infty} \frac{re^{-\beta k/n}}{1-re^{-\beta k/n}} (e^{i\theta}-1)\right)\,
&=\, \exp\left(i \sin(\theta) \sum_{k=0}^{\infty} \frac{re^{-\beta k/n}}{1-re^{-\beta k/n}} \right)\\
&\qquad \cdot
\exp\left(-2\sin^2\left(\frac{\theta}{2}\right)\sum_{k=0}^{\infty} \frac{re^{-\beta k/n}}{1-re^{-\beta k/n}}\right)
\end{split}
\end{equation*}
We will now use the result in (\ref{eq:EM2}).
Therefore, setting $\theta = \Theta/\sqrt{n}$ for $|\Theta|\leq R_n$, we have
$$
\exp\left(i \sin(\Theta/\sqrt{n}) \sum_{k=0}^{\infty} \frac{re^{-\beta k/n}}{1-re^{-\beta k/n}} \right)\,
\sim\, e^{-i n^{1/2} \beta^{-1} \ln(1-r) \Theta}\, 
=\, e^{i n^{1/2} \Theta}\, ,
$$
(note that we do not need an E-M correction because of the extra factor $1/\sqrt{n}$) and 
$$
\exp\left(-2 \sin^2\left(\frac{\Theta}{2\sqrt{n}}\right) \sum_{k=0}^{\infty} \frac{re^{-\beta k/n}}{1-re^{-\beta k/n}} \right)\,
\sim\, e^{ \beta^{-1} \ln(1-r) \Theta^2/2}\,
=\, e^{-\Theta^2/2}\, ,
$$
as $n \to \infty$, where we used $r=1-e^{-\beta}$.

So we have, for $|\Theta|<R_n$, that
\begin{equation*}
\begin{split}
\mathbf{E}\left[ 
\exp\left(\sum_{k=0}^{\infty} \frac{re^{-\beta k/n}}{1-re^{-\beta k/n}} \mathsf{T}_k (e^{i\Theta/\sqrt{n}}-1)\right)\right]\,
&=\, (1+o(1)) e^{i n^{1/2} \Theta} e^{-\Theta^2/2} \\
&\qquad
\cdot \mathbf{E}\left[ 
\exp\left(\sum_{k=0}^{\infty} \frac{re^{-\beta k/n}}{1-re^{-\beta k/n}} (\mathsf{T}_k-1) (e^{i\Theta/\sqrt{n}}-1)\right)\right]\, ,
\end{split}
\end{equation*}
as $n \to \infty$.
Now we use that $e^{i\Theta/\sqrt{n}}-1 = i\Theta \sqrt{n}(1+o(1))$ as $n \to \infty$.
Then by the Lindeberg-Feller theorem (see for example Section 3.4 of \cite{Durrett}) we have
$$
\sum_{k=0}^{\infty} \frac{re^{-\beta k/n}}{1-re^{-\beta k/n}} (\mathsf{T}_k-1) n^{-1/2}\, \stackrel{d}{\longrightarrow}\, \sigma Z\, ,
$$
as $n \to \infty$ (meaning convergence in distribution), where
$$
\sigma^2\, =\, \lim_{n \to \infty} \frac{re^{-\beta k/n}}{1-re^{-\beta k/n}} n^{-1}\, =\, \int_0^{\infty} \left(\frac{re^{-x}}{1-re^{-x}}\right)^2\, dx\, =\, e^{\beta} - 1 - \beta\, .
$$
(Exponential-1 random variables have unit variance.)
Therefore, for the characteristic function, we have
$$
\mathbf{E}\left[ 
\exp\left(\sum_{k=0}^{\infty} \frac{re^{-\beta k/n}}{1-re^{-\beta k/n}} (\mathsf{T}_k-1) (e^{i\Theta/\sqrt{n}}-1)\right)\right]\,
\to\, \mathbf{E}\left[e^{i\Theta \sigma Z}\right]\, =\, e^{-\sigma^2 \Theta^2/2}\, ,
$$
as $n \to \infty$.
Hence, putting together the two estimates we have
\begin{equation*}
\begin{split}
\mathbf{E}\left[ 
\exp\left(\sum_{k=0}^{\infty} \frac{re^{-\beta k/n}}{1-re^{-\beta k/n}} \mathsf{T}_k (e^{i\Theta/\sqrt{n}}-1)\right)\right]\,
&\sim\, e^{i\Theta \sqrt{n}} e^{-(1+\sigma^2)\Theta^2/2}\, ,
\end{split}
\end{equation*}
as $n \to \infty$.

\subsection{Completion of the calculation and consideration}

Using equation (\ref{eq:temporaryPeak}) and making the change-of-variables $\theta=\Theta/\sqrt{n}$ and $d\theta=d\Theta/\sqrt{n}$,
we have
\begin{equation*}
\begin{split}
\frac{1}{(e^{-\beta/n};e^{-\beta/n})_n}\,
	&\sim\, \frac{1}{r^n}\, \cdot \frac{1}{(r;e^{-\beta/n})_{\infty}}\, \int_{-\infty}^{\infty} e^{-in\theta}\, e^{i \Theta\, \sqrt{n}} e^{-(1+\sigma^2)\Theta^2/2}
\frac{d\Theta}{2\pi\, \sqrt{n}}\\
&=\,  \frac{1}{r^n}\, \cdot \frac{1}{(r;e^{-\beta/n})_{\infty}}\, \cdot \frac{1}{\sqrt{2\pi (1+\sigma^2) n}}\, .
\end{split}
\end{equation*}
So then using $r=1-e^{-\beta}$ and also recalling (\ref{eq:BackToIt2}), we have
\begin{equation*}
\begin{split}
\frac{1}{(e^{-\beta/n};e^{-\beta/n})_n}\,
	&\sim\, e^{n\left(-\ln(1-e^{-\beta})+\beta^{-1} \operatorname{Li}_2(1-e^{-\beta})\right)}
\, \cdot \frac{e^{-\beta/2}}{\sqrt{2\pi (1+\sigma^2) n}}\, .
\end{split}
\end{equation*}
Now in order to relate this to the Stirling formula, note that the $q$-numbers of Jackson are
$$
[n]_q\, =\, \frac{1-q^n}{1-q}\, \text{ for $n=0,1,\dots$.}
$$
So
$$
(q;q)_n\, =\, (1-q)(1-q^2)\cdots(1-q^n)\, =\, [1]_q[2]_2\cdots[n]_q \cdot (1-q)^n\, .
$$
Then writing $[n]_q! = [1]_q \cdot[2]_q \cdots [n]_q$, we have
$$
	[n]_q! \bigg|_{q=\exp(-\beta/n)}\, =\, \frac{(e^{-\beta/n};e^{-\beta/n})_n}{(1-e^{-\beta/n})^n}\, ,
$$
and $(1-e^{-\beta/n})^n \sim \left(\frac{\beta}{n}\right)^n e^{-\beta/2}$.
Thus, we have
$$
	[n]_q! \bigg|_{q=\exp(-\beta/n)}\, \sim\, n^n e^{n\left(-\ln(\beta)+\ln(1-e^{-\beta})-\beta^{-1} \operatorname{Li}_2(1-e^{-\beta})\right)}
e^{\beta}\, \sqrt{2\pi (1+\sigma^2) n}
$$
More precise asymptotics than this are available in Moak \cite{Moak}.
But this is sufficient for some purposes, such as Corollary 2.3 in \cite{StarrWalters}, noting that the integral can be performed
with the aid of a computer program such as Wolfram Mathematica.
(We wish to mention that the relation between the exponential rate and the constant prefactor does show up
in other results such as the Bahadur-Rao theorem. See for example Section 3.7 in \cite{DemboZeitouni}.)

\section{The partition function}

Hardy and Ramanujan calculated the asymptotics of the partition function, bringing some algebraic identities to bear.
Later D.~J.~E.~Newman re-derived the formula using just analysis.
Let us just apply the usual Hayman type analysis and the Euler-Maclaurin step.

The generating function is 
$$
\sum_{n=0}^{\infty} p(n) z^n\, =\, \frac{1}{\prod_{k=1}^{\infty} (1-z^k)}\, .
$$
The original Hardy-Ramanujan formula, with just the leading order asymptotics, is as follows
\begin{equation}
\label{eq:HR}
p(n)\, \sim\, \frac{1}{4n\, \sqrt{3}}\, \exp\left(\pi\, \sqrt{\frac{2n}{3}}\right)\, ,\ \text{ as  $n\to\infty$.}
\end{equation}
Let us use the Laplace transform to try to simplify the calculation.
Let $\mathsf{T}_1,\mathsf{T}_2,\dots$ be IID exponential-1 random variables. Then we can write Cauchy's integral formula as
$$
p(n)\, =\, \frac{1}{r^n}\, \mathbf{E}\left[\int_{-\pi}^{\pi} e^{-in\theta} 
\exp\left(\sum_{k=1}^{\infty} \mathsf{T}_k r^k e^{ik\theta}\right)\, \frac{d\theta}{2\pi}\right]\, ,
$$
where we have used the Fubini-Tonelli theorem.
Using the exponential-tilting formula $\mathbf{E}[e^{-a\mathsf{T_k}}f(\mathsf{T_k})] = 
(1+a)^{-1} \mathbf{E}[f((1+a)^{-1} \mathsf{T}_k)]$ for $a>-1$, we have
$$
p(n)\, =\, \frac{1}{r^n} \cdot \frac{1}{\prod_{k=1}^{\infty} (1-r^k)}\, \mathbf{E}\left[\int_{-\pi}^{\pi} e^{-in\theta} 
\exp\left(\sum_{k=1}^{\infty} \mathsf{T}_k \frac{r^k}{1-r^k}\, (e^{ik\theta}-1)\right)\, \frac{d\theta}{2\pi}\right]\, .
$$
\begin{lemma}
\label{lem:first}
Using the Euler-Maclaurin summation method, we can see that for $r=r_n=\exp(-\epsilon_n)$ with $\epsilon_n = \pi/\sqrt{6n}$,
we have
$$
\frac{1}{r^n} \cdot \frac{1}{\prod_{k=1}^{\infty} (1-r^k)}\, \sim\, 
\sqrt{\frac{\epsilon_n}{2\pi}}\, \exp\left(\pi\, \sqrt{\frac{2}{3n}}\right)\, 
=\, \frac{1}{(24n)^{1/4}}\, \exp\left(\pi\, \sqrt{\frac{2n}{3}}\right)\, ,
$$
as $n \to \infty$.
\end{lemma}
Let us delay the proof of this lemma to the end of this section. It is somewhat technical and it relies on the exact formulas for two series that happen to arise in the derivation of Stirling's formula: equations (\ref{eq:1}) and (\ref{eq:2}).
\subsection{Determination of optimal $r_n$ to leading order}
Instead of immediately proving the lemma, let us instead see why the choice of $r_n=\exp(-\epsilon_n)$ for $\epsilon_n=\pi/\sqrt{6n}$
is the optimal choice, at least to leading order.
We may approximate
$$
\prod_{k=1}^{\infty}\frac{1}{1-r^k}\, =\, \exp\left(-\sum_{k=1}^{\infty} \ln(1-r^k)\right)\, \approx\, \exp\left(-\epsilon^{-1} \int_0^{\infty} \ln(1-e^{-x})\, dx\right)\, ,
$$
where $x_k = k \epsilon$ and $\Delta x=\epsilon$ for $r=e^{-\epsilon}$. So we have
$$
\frac{1}{r^n} \cdot \frac{1}{\prod_{k=1}^{\infty} (1-r^k)}\, \approx\, 
\exp\left(\epsilon n-\epsilon^{-1} \int_0^{\infty} \ln(1-e^{-x})\, dx\right)\, ,
$$
to leading order. Of course there are lower order corrections that are captured in the Euler-Maclaurin step.
But $-\int_0^{\infty} \ln(1-e^{-x})\, dx$ is equal to $\zeta(2)=\pi^2/6$. This is one of the identities for the dilogarithm. So this gives $\epsilon = \pi / \sqrt{6n}$, and to
leading order this then gives
$$
\frac{1}{r^n} \cdot \frac{1}{\prod_{k=1}^{\infty} (1-r^k)}\, \approx\, 
\exp\left(\pi\, \sqrt{\frac{2n}{3}}\right)\, .
$$
We will prove the full Lemma \ref{lem:first}, later.

\subsection{Treatment of main contribution: major arc}
Let us add-and-subtract the expectation of each random variable to obtain
$$
p(n)\, =\, \frac{1}{r^n} \cdot \frac{1}{\prod_{k=1}^{\infty} (1-r^k)}\, \mathbf{E}\left[\int_{-\pi}^{\pi} e^{-in\theta+\sum_{k=1}^{\infty} \frac{r^k}{1-r^k}\, (e^{ik\theta}-1)}
\exp\left(\sum_{k=1}^{\infty} (\mathsf{T}_k-1) \frac{r^k}{1-r^k}\, (e^{ik\theta}-1)\right)\, \frac{d\theta}{2\pi}\right]\, .
$$
Now the relevant integrals are
$$
	\int_0^{\infty} \frac{x}{e^x-1}\, dx\, =\, \frac{\pi^2}{6}\, ,
$$
(which is related by IBP to $-\int_0^{\infty} \ln(1-e^{-x})\, dx$)
and
$$
	\int_0^{\infty} \frac{x^2}{e^x-1}\, dx\, =\, 2\zeta(3)\, .
$$
Using that, and letting $\eta=\eta_n=\epsilon_n^{3/2}/\sqrt{2\zeta(3)}=\frac{2\pi^{3/2}}{\sqrt{\zeta(3)}\, (24n)^{3/4}}$, we have
$$
	\eta \sum_{k=1}^{\infty} \frac{r^k}{1-r^k}\, k\, =\, \eta \epsilon^{-2} \sum_{k=1}^{\infty} \frac{\epsilon k}{e^{\epsilon k}-1}\, \epsilon\,
	\sim\, \eta \epsilon^{-2}\, \frac{\pi^2}{6}\, \sim\, \eta n\, ,
$$
and
$$
	\eta^2 \sum_{k=1}^{\infty} \frac{r^k}{1-r^k}\, k^2\, 
	=\, \eta^2 \epsilon^{-3} \sum_{k=1}^{\infty} \frac{(\epsilon k)^2}{e^{\epsilon k}-1}\, \epsilon
	\sim\, 2\zeta(3) \eta^2 \epsilon^{-3}\, =\, 1\, .
$$
Therefore, setting $\theta = \Theta \eta$ for $-R_n \leq \Theta \leq R_n$, for $R_n \to \infty$ as $n \to \infty$ {\em slowly}, we have
$$
\sum_{k=1}^{\infty} \frac{r^k}{1-r^k}\, (i \sin(k\theta) - 2 \sin^2(k\theta/2))\,
	=\, \eta n \Theta - \frac{\Theta^2}{2} + o(1)\, ,
$$
as $n \to \infty$. One should do the Euler-Maclaurin step to check that there is no extra order-1 or larger correction. But there is not.
So we have
$$
p(n)\, \sim\, \frac{1}{r^n} \cdot \frac{1}{\prod_{k=1}^{\infty} (1-r^k)}\, 
\left(\mathbf{E}\left[\int_{-R_n}^{R_n} e^{-\Theta^2/2}
\exp\left(\sum_{k=1}^{\infty} (\mathsf{T}_k-1) \frac{r^k}{1-r^k}\, (e^{ik\eta\Theta}-1)\right)\, \frac{d\Theta}{2\pi/\eta}\right]
+ \mathcal{E}_n\right)\, ,
$$
where $\mathcal{E}_n$ is the remainder of the integral for $|\theta| \in ( \eta R_n,\pi)$. 
Finally, we note that
$$
	\operatorname{Var}\left(\sum_{k=1}^{\infty} (\mathsf{T}_k-1) \frac{\epsilon k}{r^{-k}-1}\, \epsilon^{1/2}\right)\, 
	\sim\, \int_0^{\infty} \left(\frac{x}{e^x-1}\right)^2\, dx\, =\, \frac{\pi^2}{3} - 2\zeta(3)\, . 
$$
Since the random variable is well-approximated by a normal random variable, the CLT and the characteristic equation formula
for a normal random variable give
\begin{equation*}
\begin{split}
\mathbf{E}\left[\exp\left(\sum_{k=1}^{\infty} (\mathsf{T}_k-1) \frac{r^k}{1-r^k}\, ik\eta\Theta\right)\right]\,
&\sim\, \exp\left(-\eta^2\epsilon^{-3}\left(\frac{\pi^2}{3}-2\zeta(3)\right)\Theta^2/2\right) \\
&=\, \exp\left(-\frac{\pi^2}{12 \zeta(3)} \Theta^2\right) e^{\Theta^2/2}\, .
\end{split}
\end{equation*}
So we obtain
$$
p(n)\, \sim\, \frac{1}{r^n} \cdot \frac{1}{\prod_{k=1}^{\infty} (1-r^k)}\, 
\left(\int_{-R_n}^{R_n} \exp\left(-\frac{\pi^2}{12 \zeta(3)} \Theta^2\right)\, \frac{d\Theta}{2\pi/\eta}
+ \mathcal{E}_n\right)\, ,
$$
But this is 
$$
p(n)\, \sim\, \frac{1}{r^n} \cdot \frac{1}{\prod_{k=1}^{\infty} (1-r^k)}\, 
\left(\frac{\eta}{2\pi}\, \sqrt{2\pi\, \frac{6\zeta(3)}{\pi^2}}
+ \mathcal{E}_n\right)\, ,
$$
But since $\eta=\frac{2\pi^{3/2}}{\sqrt{\zeta(3)}\, (24n)^{3/4}}$, this gives
$$
p(n)\, \sim\, \frac{1}{r^n} \cdot \frac{1}{\prod_{k=1}^{\infty} (1-r^k)}\, 
\left(2\, \sqrt{\frac{\pi^3}{\zeta(3)\cdot (24n)^{3/2}}}\cdot\, \sqrt{\frac{3\zeta(3)}{\pi^3}}
+ \mathcal{E}_n\right)\, ,
$$
as $n \to \infty$.
Therefore,
$$
p(n)\, \sim\, \frac{1}{r^n} \cdot \frac{1}{\prod_{k=1}^{\infty} (1-r^k)}\, 
\left(\frac{2\, \sqrt{3}}{(24n)^{3/4}}+ \mathcal{E}_n\right)\, ,
$$
as $n \to \infty$.
Using Lemma \ref{lem:first} this establishes the result in equation (\ref{eq:HR}), as long as we prove that
$\mathcal{E}_n = o(n^{-3/4})$.

\subsection{Bounding corrections: minor arc and the lemma}

\begin{prop}[Cram\'er's theorem] The sequence $(\mathsf{T}_1+\dots+\mathsf{T}_n)/n$ satisfies a large deviation principle
with large deviation rate function $\Lambda^*$, whose formula is 
$$
\Lambda^*(x)\, =\, x-1-\ln(x)\, ,
$$
the Legendre transform of the  cumulant generating function $\Lambda(t) = \ln\left(\mathbf{E}[e^{t\mathsf{T}}]\right)=-\ln(1-t)$.
\end{prop}

Using that, we can see that, for $t = o(\epsilon^{-1})$, we have
$$
	2\sum_{k=1}^{\infty} \mathsf{T}_k\, \frac{r^k}{1-r^k}\, \sin^2\left(\frac{t \varepsilon k}{2}\right)\,
	\sim\, 2\epsilon^{-1} \int_0^{\infty} \frac{1}{e^x-1} \sin^2\left(\frac{t x}{2}\right)\, dx\, ,
$$
with exponentially small errors controlled by the large deviation rate function.
Moreover, for $t >\epsilon^{-3/4}$ taking $\theta = \epsilon t >\epsilon^{1/4}$ also satisfying $|\theta|<\pi$, we have 
$2\mathsf{T}_1 \frac{r}{1-r} \sin^2(\theta/2)\geq 2\varepsilon^{-1/2}\mathsf{T}_1/\pi^2$ times $(1+o(1))$,
as $n \to \infty$. Since
$$
	\left|e^{\sum_{k=1}^{\infty} \mathsf{T}_k\, \frac{r^k}{1-r^k} (e^{i\theta}-1)}\right|\, =\, 
	\exp\left(-2\sum_{k=1}^{\infty} \mathsf{T}_k\, \frac{r^k}{1-r^k}\, \sin^2(\theta k/2)\right)\, ,
$$
this establishes that $\mathcal{E}_n$ equals $\exp(-Cn^{1/4})$ plus
$$
2\int_{R_n \eta/\epsilon}^{\epsilon^{-3/4}} 
\exp\left(- 2\epsilon^{-1} \int_0^{\infty} \frac{1}{e^x-1} \sin^2\left(\frac{t x}{2}\right)\, dx\right)\, 
\frac{dt}{\epsilon^{-1}}\, =\, O\left(\epsilon e^{-CR_n^2}\right)\, .
$$
So as long as $CR_n^2/\ln(n) >1/4$, we have the desired result that $\mathcal{E}_n = o(n^{-3/4})$.
So we choose $R_n = \sqrt{\ln(n)/C}$.

\begin{proofof}{\bf Proof of Lemma \ref{lem:first}:}
We note that
$$
\ln\left(\frac{1}{\prod_{k=1}^{\infty}(1-r^k)}\right)\, =\, -\sum_{k=1}^{\infty} \ln(1-r^k)\, .
$$
Then by the Euler-Maclaurin method, we have
\begin{equation*}
\begin{split}
	-\sum_{k=1}^{\infty} \ln(1-r^k)\, &=\, 
	-\sum_{k=1}^{\infty} \ln(1-e^{-\epsilon k})\\
	&=\, -\int_{1/2}^{\infty} \ln\left(1-e^{-\epsilon x}\right)\, dx
	+\sum_{k=1}^{\infty} \int_{-1/2}^{1/2} \ln\left(\frac{1-e^{-\epsilon k}e^{-\epsilon x}}{1-e^{-\epsilon k}}\right)\, dx\, .
\end{split}
\end{equation*}
Letting $y=\epsilon x$, so that $dx = dy/\epsilon$, we see that
\begin{equation*}
-\int_{1/2}^{\infty} \ln\left(1-e^{-\epsilon x}\right)\, dx\,
=\, -\epsilon^{-1} \int_0^{\infty} \ln(1-e^{-y})\, dy + \epsilon^{-1} \int_0^{\epsilon/2} \ln(1-e^{-y})\, dy\, .
\end{equation*}
But since $-\int_0^{\infty} \ln(1-e^{-y})\, dy$ equals $\zeta(2)=\pi^2/6$ we may then rewrite this as 
\begin{equation*}
-\int_{1/2}^{\infty} \ln\left(1-e^{-\epsilon x}\right)\, dx\,
=\, -\epsilon^{-1}\, \frac{\pi^2}{6} + \epsilon^{-1} \int_0^{\epsilon/2} \ln(1-e^{-y})\, dy\, .
\end{equation*}
Then we may expand
$$
\epsilon^{-1} \int_0^{\epsilon/2} \ln(1-e^{-y})\, dy\,
=\, \epsilon^{-1} \int_0^{\epsilon/2} \ln(y)\, dy + o(1)\,
=\, \frac{1}{2}\, \ln\left(\frac{\epsilon}{2}\right) - \frac{1}{2} +o(1)\, .
$$
So we have
\begin{equation*}
\begin{split}
	-\sum_{k=1}^{\infty} \ln(1-e^{-\epsilon k})\,
	&=\, \epsilon^{-1}\, \frac{\pi^2}{6}
+ \frac{1}{2}\, \ln\left(\frac{\epsilon}{2}\right) - \frac{1}{2} +o(1)
	+\sum_{k=1}^{\infty} \int_{-1/2}^{1/2} \ln\left(\frac{1-e^{-\epsilon k}e^{-\epsilon x}}{1-e^{-\epsilon k}}\right)\, dx\, .
\end{split}
\end{equation*}
Now let us consider the series. We may rewrite
\begin{equation*}
\begin{split}
	\sum_{k=1}^{\infty} \int_{-1/2}^{1/2} \ln\left(\frac{1-e^{-\epsilon k}e^{-\epsilon x}}{1-e^{-\epsilon k}}\right)\, dx\, 
	&=\, 	\sum_{k=1}^{\infty} \int_{-1/2}^{1/2} \ln\left(1+\frac{e^{-\epsilon k}(1-e^{-\epsilon x})}{1-e^{-\epsilon k}}\right)\, dx\\
	&=\, 	\sum_{k=1}^{\infty} \int_{-1/2}^{1/2} \ln\left(1+\frac{e^{-\epsilon k}\epsilon x}{1-e^{-\epsilon k}}\right)\, dx\\
&\qquad + \sum_{k=1}^{\infty} \int_{-1/2}^{1/2} \ln\left(1-\frac{1}{2}\, \cdot \frac{e^{-\epsilon k}(\epsilon x)^2}
{1-e^{-\epsilon k}+e^{-k\epsilon} \epsilon x}\right)\, dx + o(1)\, .
\end{split}
\end{equation*}
For the first series, we note that
\begin{equation*}
\begin{split}
\int_{-1/2}^{1/2} \ln\left(1+\frac{e^{-\epsilon k}\epsilon x}{1-e^{-\epsilon k}}\right)\, dx\,
&\sim\, \int_{-1/2}^{1/2} \ln\left(1+\frac{x}{k}\right)\, dx\\
&=\, \ln\left(1-\frac{1}{4k^2}\right) + \left(-1+k\, \ln\left(1+\frac{1}{2k}\right)-k\, \ln\left(1-\frac{1}{2k}\right)\right)\, .
\end{split}
\end{equation*}
Using equations (\ref{eq:1}) and (\ref{eq:2}), this gives
$$
\sum_{k=1}^{\infty} \int_{-1/2}^{1/2} \ln\left(1+\frac{e^{-\epsilon k}\epsilon x}{1-e^{-\epsilon k}}\right)\, dx\,
\sim\, \frac{1-\ln(\pi)}{2}\, .
$$
Thus, to establish the lemma, we just need to show that the second series converges to $0$.

For the second series, we have
\begin{equation*}
\begin{split}
	\int_{-1/2}^{1/2} \ln\left(1-\frac{1}{2}\, \cdot \frac{\epsilon^2 x^2}
	{e^{\epsilon k}-1+ \epsilon x}\right)\, dx\,
	&\sim\, -\frac{1}{2}\, \int_{-1/2}^{1/2} \frac{\epsilon^2 x^2}{e^{\epsilon k}-1+\epsilon x}\, dx\, \sim\, -\frac{\epsilon^2}{24(e^{\epsilon k}-1)}\, .
\end{split}
\end{equation*}
But then that means
$$
\sum_{k=1}^{\infty} \int_{-1/2}^{1/2} \ln\left(1-\frac{1}{2}\, \cdot \frac{e^{-\epsilon k}(\epsilon x)^2}
{1-e^{-\epsilon k}+e^{-k\epsilon} \epsilon x}\right)\, dx\,
\sim\, -\frac{\epsilon}{24} \int_{\epsilon}^{\infty} \frac{1}{e^{x}-1}\, dx\, .
$$
But this is $O(\epsilon \ln(\epsilon))$. In other words, it converges to $0$ as $n \to \infty$.
\end{proofof}

\subsection{Completion of the argument and consideration}

Using this with (\ref{eq:temporaryPeak}), we have, making the change-of-variables $\theta=\Theta/\sqrt{n}$ so that $d\theta=d\Theta/\sqrt{n}$,

\appendix

\section{Technical lemmas for Newton's binomial identity}
\label{app:first}

\begin{proofof}{\bf Proof of Lemma \ref{lem:Hay1}:}
First note that $\exp(-2n\tau \sin^2(\theta/2)) \leq \exp(-2n \tau \theta^2/\pi^2)$.
So let $R_n = \pi \sqrt{\ln(n)}$ and note that for $|\theta|>R_n/\sqrt{n\tau}$, we have 
$\exp(-2n\tau \sin^2(\theta/2)) \leq 1/n^2$.
Therefore, setting $\Theta = \theta \sqrt{n\tau}$, so that $d\theta=d\Theta/\sqrt{n\tau}$, we have
$$
	\int_{-\pi}^{\pi} \exp\left(-in\theta-n\tau (e^{i\theta}-1)\right)\, \frac{d\theta}{2\pi}\,
	=\, \int_{-R_n}^{R_n} e^{-i\Phi_n(\tau,\Theta)}\, e^{-A_n(\tau,\Theta)}\, \frac{d\Theta}{2\pi\, \sqrt{n\tau}}
+ O\left(\frac{1}{n^2}\right)\, ,
$$
where
$$
A_n(\tau,\Theta) = 2n\tau \sin^2\left(n^{-1/2} \tau^{-1/2} \Theta/2\right)\, ,
$$
and
$$
\Phi_n(\tau,\Theta)\, =\, n^{1/2} \tau^{-1/2} \left(\Theta- n^{1/2} \tau^{3/2} \sin\left(\Theta n^{-1/2} \tau^{-1/2}\right)\right)\, .
$$
Now consider the open region of $\C$
$$
	\mathcal{R}_n(\tau)\, =\, \{\Theta+i\rho n^{1/2} \tau^{1/2} \ln(\tau)\, :\, -R_n<\Theta<R_n\, ,\ 0<\rho<1\}\, ,
$$
which is a bounded rectangle.
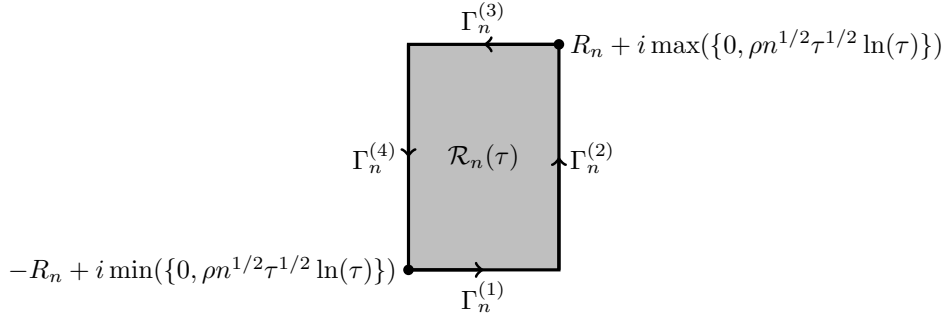
\begin{figure}
\begin{center}
\begin{tikzpicture}[very thick]
	\fill[black!25!white] (0,0) rectangle (2,3);
	\draw (1,1.5) node[] {\small $\mathcal{R}_n(\tau)$};
	\draw[->] (0,0) -- (1,0)  node[below] {\small $\Gamma_n^{(1)}$};
	\draw[->] (2,0) -- (2,1.5) node[right] {\small $\Gamma_n^{(2)}$};
	\draw[->] (2,3) -- (1,3) node[above] {\small $\Gamma_n^{(3)}$};
	\draw[->] (0,3) -- (0,1.5)  node[left] {\small $\Gamma_n^{(4)}$};
	\draw (0,0) rectangle (2,3);
	\fill (0,0) circle (2pt) node[left] {\small $-R_n+i\min(\{0,\rho n^{1/2} \tau^{1/2} \ln(\tau)\})$};
	\fill (2,3) circle (2pt) node[right] {\small $R_n+i\max(\{0,\rho n^{1/2} \tau^{1/2} \ln(\tau)\})$};
\end{tikzpicture}
\end{center}
\caption{A schematic of the contour $\Gamma_n(\tau)$ for the application of Cauchy's integral formula and the region $\mathcal{R}_n(\tau)$.
\label{fig:contour}}
\end{figure}
Let $\Gamma_n(\tau)$ be the counterclockwise oriented contour around the boundary of $\mathcal{R}_n(\tau)$.
Let the line segments be $\Gamma_n^{(1)}$, $\Gamma_n^{(2)}$, $\Gamma_n^{(3)}$ and $\Gamma_n^{(4)}$ as shown in Figure \ref{fig:contour}.
Note that if $\tau < 1$ then we will obtain the opposite of the integral we want in $\Gamma_n^{(3)}$. But that is not the point:
we have 
$$
	\oint_{\Gamma_n(\tau)} \exp\left(-in^{1/2} \tau^{-1/2} z + n \tau (e^{izn^{-1/2}\tau^{-1/2}}-1)\right)\, dz\, =\, 0\, ,
$$
by the Cauchy integral formula. (See for example \cite{MarsdenHoffman}.)
So our immediate goal is just to show that the integral over $\Gamma_n^{(1)}$ and $\Gamma_n^{(3)}$ are nearly opposite with
small corrections. In other words, let us first show that the integrals over $\Gamma_n^{(2)}$ and $\Gamma_n^{(4)}$ are small.

A calculation using Taylor series shows that
\begin{equation*}
\begin{split}
z=\Theta+i\rho n^{1/2}\tau^{1/2} \ln(\tau)\qquad \Rightarrow\qquad\\
&\hspace{-4.5cm}
-in^{1/2} \tau^{-1/2} z + n \tau (e^{izn^{-1/2}\tau^{-1/2}}-1)\, =\, \mathcal{X}_n(\tau,\rho,\Theta) + o(1)\, ,
\end{split}
\end{equation*}
where
$$
\mathcal{X}_n(\tau,\rho,\theta)\, =\, n\rho \ln(\tau) + n\tau(\tau^{-\rho}-1) + i n^{1/2} (\tau^{\frac{1}{2}-\rho}-\tau^{-1/2})\Theta
- \frac{1}{2}\, \tau^{-\rho} \Theta^2\, .
$$
(We used $\tau^{-1} < n^{\delta}$ with $\delta<1/3$.)
We note that $ \mathcal{X}_n(\tau,1,\Theta)$ is a pure quadratic in $\Theta$. This means that the integrand in either $\Gamma_n^{(1)}$ or $\Gamma_n^{(3)}$ is a pure Gaussian.  
That motivates the choice of this contour $\Gamma_n(\tau)$.
In other words, if we can bound the integrals over $\Gamma_n^{(2)}$ and $\Gamma_n^{(4)}$ then we can show that
$$
\int_{-R_n}^{R_n} e^{-i\Phi_n(\tau,\Theta)}\, e^{-A_n(\tau,\Theta)}\, \frac{d\Theta}{2\pi\, \sqrt{n\tau}}\, 
=\, (1+o(1))e^{n \left(\ln(\tau)+1-\tau\right)} \int_{-R_n}^{R_n} e^{-\Theta^2/(2\tau)}\, \frac{d\Theta}{2\pi\, \sqrt{n\tau}} + \mathcal{E}_n\, ,
$$
where $\mathcal{E}_n$ is the error arising from bounding the integrals along $\Gamma_n^{(2)}$ and $\Gamma_n^{(4)}$.

For $0<\rho<1$ we have
\begin{equation*}
\begin{split}
z\, =\, i\rho n^{1/2}\tau\ln(\tau) \pm R_n\qquad \Rightarrow\\
&\hspace{-3cm}
\left|e^{-in^{1/2} \tau^{-1/2} z + n \tau (e^{izn^{-1/2}\tau^{-1/2}}-1)}\right|\,
\sim\, e^{n\rho \ln(\tau) + n \tau(\tau^{-\rho}-1)-\frac{1}{2}\, \tau^{-\rho} R_n^2}\, .
\end{split}
\end{equation*}
By convexity, 
$$
n\rho \ln(\tau) + n \tau(\tau^{-\rho}-1)\, \leq\, n\rho (\ln(\tau)+1-\tau)\, .
$$
If $\tau<1$ then $\tau^{-\rho} R_n^2 \geq R_n^2$.
So, if $\tau<1$ then
\begin{equation*}
\begin{split}
\left|\int_{\Gamma^{(2,4)}_n}  \exp\left(-in^{1/2} \tau^{-1/2} z + n \tau (e^{izn^{-1/2}\tau^{-1/2}}-1)\right)\, \frac{dz}{2\pi \sqrt{n\tau}}\right|\\
&\hspace{-3cm}
\leq\, \left(1+o(1)\right) \frac{|\ln(\tau)|}{2\pi}\, e^{-R_n^2/2} \int_0^1e^{n\rho (\ln(\tau)+1-\tau)}\, d\rho\, .
\end{split}
\end{equation*}
If $\tau<1/e$ then this is bounded by $(1+o(1))e^{-R_n^2/2}e/(2\pi)$ which is $O(1/n^2)$ (by continuing the domain of integration to $(0,\infty)$ and doing the exponential integral).
If $1/e<\tau<1$ then we can bound the integrand by just $1$ and we obtain a bound of $(1+o(1))e^{-R_n^2/2}/(2\pi)$ which is $O(1/n^2)$.

Now let us consider $\tau>1$. If $|\tau-1|$ is small, then we essentially are in the same situation as when we had $1/e<\tau<1$.
For $\tau\gg 1$, let us use
\begin{equation*}
\begin{split}
\left|\int_{\Gamma^{(2,4)}_n}  \exp\left(-in^{1/2} \tau^{-1/2} z + n \tau (e^{izn^{-1/2}\tau^{-1/2}}-1)\right)\, \frac{dz}{2\pi \sqrt{n\tau}}\right|\\
&\hspace{-3cm}
\leq\, \left(1+o(1)\right) \frac{|\ln(\tau)|}{2\pi}\, \int_0^1e^{n\rho (\ln(\tau)+1-\tau)}\, d\rho\, .
\end{split}
\end{equation*}
This is bounded by $O(1) |\ln(\tau)|/(n(\tau-1-\ln(\tau)))$ which is $o(1/n)$(by continuing the domain of integration to $(0,\infty)$ and doing the exponential integral). 
That is sufficient for our purposes.

So we have proved that in all cases we have
$$
\int_{-R_n}^{R_n} e^{-i\Phi_n(\tau,\Theta)}\, e^{-A_n(\tau,\Theta)}\, \frac{d\Theta}{2\pi\, \sqrt{n\tau}}\, 
=\, (1+o(1))e^{n \left(\ln(\tau)+1-\tau\right)} \int_{-R_n}^{R_n} e^{-\Theta^2/(2\tau)}\, \frac{d\Theta}{2\pi\, \sqrt{n\tau}} + o\left(\frac{1}{n}\right)\, .
$$
Then by a simple change of variables,
$$
\int_{-R_n}^{R_n} e^{-\Theta^2/(2\tau)}\, \frac{d\Theta}{2\pi\, \sqrt{n\tau}}\,
	=\, \int_{-R_n/\sqrt{\tau}}^{R_n/\sqrt{\tau}} e^{-x^2/2}\, \frac{dx}{2\pi\, \sqrt{n}}\, .
$$
If $\tau\gg 1$ so that the domain of integration is much smaller than $(-R_n,R_n)$, then we also have $e^{n \left(\ln(\tau)+1-\tau\right)}\ll e^{-n}$ which is
$o(1/n)$.
This proves the lemma.
\end{proofof}

\begin{proofof}{\bf Proof of Lemma \ref{lem:Hay2}:}
We note that, by Taylor expansion
$$
\ln(\tau) + 1 - \tau\, =\, -\frac{(\tau-1)^2}{2} + O(|\tau-1|^3)\, ,\ \text{ as $\tau \to 1$.}
$$
Moreover it is concave down. Since $e^{-s\tau} \tau^{s-1}$ is integrable, it suffices to note that for a suitable cutoff $R_n'$ we may neglect $e^{n(\ln(\tau)+1-\tau)}$
as being $o(1/n)$ for $|\tau-1|n^{-1/2}>R_n'$. And this quantity is $O(\sqrt{\ln(n)})$, where for $|\tau-1|n^{-1/2}<R_n'$ we may just approximate by $e^{n(\ln(\tau) + 1 - \tau)}
\sim e^{n(\tau-1)^2/2}$.
\end{proofof}

\begin{proofof}{Proof of Corollary \ref{cor:Gamma}:}
We know that for $\tau < n^{-\delta}$ the integrand is $O(n^{-\epsilon n})$. But elsewhere we have the integral is asymptotic to $\frac{1}{n}\, (1+o(1))$. So the result
follows because $O(n^{-\epsilon n})$ is clearly $o(1/n)$.
\end{proofof}

\section{Summary and outlook}

In the process of determining moments of collisions of 1d random walks, the author along with Samen Hossein was led to the use of the Laplace
transform in the context of extracting asymptotics from generating functions.
The usual approach of Hayman's method, approximating the integrand in Cauchy's integral formula by a Gaussian, is easiest
to apply if the generating function is an exponential type function. (For example, extracting the asymptotics of $1/n!$ from the generating function
$\exp(z)$ is the most basic model.) If the generating function is closer to a rational function, for example a finite product of inverse-linear factors, or a convergent infinite product of such factors, then it seems useful to apply the Laplace transform method.

Without mentioning the Laplace transform, this approach was already carried out by Fristedt and later Romik, using infinitely many
exponential-1 random variables. See \cite{Fristedt} and \cite{Romik}
We note that the results of this paper are essentially an exercise in Romik's textbook \cite{RomikText}.
But an important element of his suggestion is to use D.~J.~Newman's rewriting of the generating function, which
we do not find to be the most naive approach \cite{Newman}.
Similarly, one should refer to the original paper of Ramanujan and Hardy \cite{HardyRamanujan}.
But we find the use of the pentagonal identity to be somewhat mysterious for newcomers, as opposed to simply applying 
the Euler-Maclaurin step and Hayman's method (as is used in more elementary asymptotic analyses).
Similarly, Erd\"os has a purely analytical argument, since he eschewed complex analysis \cite{Erdos}.
But the complex-analytic approach is the most common approach, such as what most students see to derive Stirling's formula.
Newman has done just what we would want.
But it seems to go more smoothly by using exponential-1 random variables, following Fristedt and Romik.

Of course Hardy and Ramanujan's approach is justified by their full asymptotic expansion, later improved by
Rademacher into a convergent series \cite{Rademacher}. But here we just seek the leading order
asymptotics.

For our approach, all that is needed is the knowledge of complex analysis and a first course in graduate probability, roughly commensurate
with a first year graduate student's knowledge or beyond.
One does need to do various integrals, especially those that lead to the dilogarithm function.
For this, it helps to use a program such as Wolfram's Mathematica.
We note that Ramanujan, himself, proved some of the known dilogarithm identities. 
See for example Berndt \cite{Berndt}.

Two other examples of the Laplace transform's utility in extracting asymptotics were provided applied to the generating functions: Newton's 
binomial formula and the $q$-binomial formula for $a=0$. These are easier, and were provided as warm-up examples.

\section*{Acknowledgments}
I am grateful to Timothy Li who contributed many useful discussions.

\label{sec:Outlook}

\baselineskip=12pt
\bibliographystyle{plain}

\end{document}